\documentclass[10pt,oneside,a4paper]{amsart}
\usepackage{amsmath,amsfonts,amsthm, amssymb}
\usepackage{array}
\usepackage[all,textures]{xy}

\theoremstyle{plain}
\newtheorem{theorem}{Theorem}[section]
\newtheorem{proposition}[theorem]{Proposition}
\newtheorem{lemma}[theorem]{Lemma}
\newtheorem{corollary}[theorem]{Corollary}

\theoremstyle{definition}
\newtheorem{definition}[theorem]{Definition}

\theoremstyle{remark}

\newtheorem{remark}[theorem]{Remark}
\newtheorem{example}[theorem]{Example}

\newcommand{\ko}{\: , \;}

\newcommand{\ie}{{\em i.e.\ }}

\newcommand{\id}{\mathbf{1}}

\newcommand{\Cokern}{\mathrm{Coker}}

\renewcommand{\lim}{\mathrm{lim}}

\newcommand{\Hom}{\mathrm{Hom}}

\newcommand{\Z}{\mathbb{Z}}

\newcommand{\N}{\mathbb{N}}

\newcommand{\CC}{\mathbf{C}}

\newcommand{\Mod}{\ensuremath{\mathsf{Mod}} }
\newcommand{\Com}{\ensuremath{\mathsf{Com}} }
\newcommand{\Gr}{\ensuremath{\mathsf{Gr}} }

\newcommand{\lra}{\longrightarrow}
\newcommand{\ra}{\rightarrow}

\newcommand{\aaa}{\ensuremath{\mathcal{A}}}

\newcommand{\mmm}{\ensuremath{\mathcal{M}}}

\newcommand{\ppp}{\ensuremath{\mathcal{P}}}

\CompileMatrices
\SelectTips{cm}{11}

\title[Vanishing derived categories]{On the (non)vanishing of
some ``derived'' categories of curved dg algebras}
\author{Bernhard Keller}
\address[Bernhard Keller]{Universit\'e Paris Diderot - Paris 7,
Institut de Math\'ematiques de Jussieu,
U.M.R. 7586 du CNRS,
UFR de Math\'ematiques, Case 7012,
B\^{a}timent Chevaleret,
75205 Paris Cedex 13,
FRANCE}
\email{keller@math.jussieu.fr}

\author{Wendy Lowen$^{1}$}
\address[Wendy Lowen]{Departement Wiskunde-Informatica, Middelheimcampus,
Middelheimlaan 1,
2020 Antwerp, Belgium}
\email{wendy.lowen@ua.ac.be}

\author{Pedro Nicol\'as$^{2}$}
\address[Pedro Nicol\'as]{Institut de Math\'ematiques de Jussieu,
U.M.R. 7586 du CNRS,
Universit\'e Paris Diderot - Paris 7,
UFR de Math\'ematiques, Case 7012,
B\^{a}timent Chevaleret,
75205 Paris Cedex 13,
FRANCE}
\email{nicolas@math.jussieu.fr}

\thanks{$^{1}$Postdoctoral Fellow of the Research Foundation Flanders (FWO). $^{2}$Postdoctoral Fellow of the Fundaci\'on S\'eneca de la Regi\'on de Murcia.}

\begin{document}
\maketitle

\begin{abstract}
  Since curved dg algebras, and modules over them, have differentials
  whose square is not zero, these objects have no cohomology, and there
  is no classical derived category. For different purposes, different
  notions of ``derived'' categories have been introduced in the
  literature. In this note, we show that for some concrete curved dg
  algebras, these derived categories vanish. This happens for example
  for the initial curved dg algebra whose module category is the
  category of precomplexes, and for certain deformations of dg
  algebras.
\end{abstract}

\section{Introduction}

Curved dg algebras and modules were introduced in \cite{positselski}, in relation with quadratic duality. Examples of a different nature occur as deformations of ordinary dg algebras. Indeed, inspection of the Hochschild complex of a dg algebra immediately reveals the possible occurrence of curvature in deformations. The deformation theory of algebras (\cite{gerstenhaber, gerstenhaber1}) and of abelian categories (\cite{lowenvandenbergh1, lowenvandenbergh2}) suggests that deformation should somehow take place on the derived level.
For derived categories of abelian categories, the situation was investigated in \cite{lowen7}. For dg algebras, the relation between Hochschild cohomology and derived Morita deformations was investigated in \cite{kellerlowen}, where it was shown that not every Hochschild cocycle can be realized by means of a Morita deformation of the dg algebra. This raises further questions as to the possibility of deriving deformed curved dg algebras. More precisely: suppose $\bar{A}$ is a curved dg algebra deforming an ordinary dg algebra $A$, is there a reasonable definition of derived category $D_?(\bar{A})$ which can be considered to be a ``derived deformation'' of $D(A)$? Since curved dg algebras fail to have square zero differentials, and hence fail to have cohomology objects, a straightforward generalization of the definition of the derived category of a dg algebra does not exist. Different candidate derived categories have been considered in the literature (\cite{nicolas}, \cite{po!
 sitselski2}), but none of these is such that for all dg algebras, the newly defined category coincides with the classical derived category.

Our answer to the general existence of ``derived deformations'' is a negative one: we give examples where it is impossible to define a reasonable derived category $D_?(\bar{A})$ deforming $D(A)$. By \emph{reasonable}, we mean satisfying some combination of a number of natural axiomatic requirements (listed in \ref{parreason}) for the corresponding class of ``acyclic'' objects. Loosely speaking, we will refer to these categories as ``derived'' categories. By \emph{deforming}, we mean that a complex over $A$ is acyclic if and only if its image over $\bar{A}$ is ``acyclic''.   Our most pronounced example in this respect is the ``graded field'' $A = k[u, u^{-1}]$ where $u$ is of degree $2$. The element $u$ gives rise to a Hochschild cocycle and an infinitesimal deformation $\bar{A}$, but there is no ``derived'' category $D_?(\bar{A})$ deforming $D(A)$. Moreover, over a field, the only ``derived'' category of $\bar{A}$ is actually zero (Proposition \ref{prophor}). Another class o!
 f curved dg algebras whose ``derived'' categories we show to vanish, are the ``initial cdg algebras'' $k[c]$ and $k[c]/c^n$ for $c$ in degree $2$ and $n \geq 2$ (Proposition \ref{derivedzero}). In \S \ref{parmainhop}, we take a slightly different approach to ``derived'' categories, by looking at classes of ``homotopical projectives''. The existence of non-zero ``derived'' categories is closely related to the existence of graded projective, respectively graded projective and graded small objects in the homotopy category. In particular, we show that the deformation $\bar{A}$ of $A = k[u]$ (with $u$ in degree $2$) corresponding to the cocycle $u$, posesses a non-zero ``derived'' category $D_?(\bar{A})$, but this category actually does not deform the classical $D(A)$ (Proposition \ref{propnondef}).

Finally, in \S \ref{parlit}, we take a closer look at particular candidate derived categories studied in the literature. In \S \ref{parbar}, we look at the bar derived category $D_{\mathrm{bar}}(A)$ of \cite{nicolas}, which is defined for a unital cdg algebra $A$ over an arbitrary commutative ground ring $k$,  and which should be regarded as a curved analogue of the relative derived category of a dg algebra (in which by definition the $k$-contractible complexes become zero). We show that if $k$ is a field and $A$ has a non-zero curvature, $D_{\mathrm{bar}}(A) = 0$. This is a consequence of the fact that the bar derived category vanishes for the initial cdg algebras, and satisfies a strong base change property (see \S \ref{parbase}). In \S \ref{parsecond}, we take a look at the ``derived categories of the second kind'' defined in \cite{positselski2}. These categories (of which there are three subtypes) can be regarded as universal ``derived'' categories. The existence of non-!
 zero derived categories of the second kind over a field, in spite of their vanishing on $k[c]$, can be explained by the fact that they don't satisfy the strong base change property.

In contrast to the approaches in \cite{nicolas} and \cite{positselski2}, which make use of the interplay between algebras and coalgebras through the bar/cobar formalism, the methods in this paper are elementary (except in \S \ref{parbar} where we apply our results to the setting of \cite{nicolas}).

\subsection{Acknowledgement}
Most of this work, in particular the vanishing of a number of
``derived categories'' satisfying certain natural axioms (and hence,
of the bar derived category over a field) originated in 2006, when
the authors were together in Paris. However, it was not until we
discovered the beautiful applications of derived categories
satisfying precisely those axioms in the work of Leonid Positselski,
as presented by him in Paris in April 2009, that we decided it would
be worthwile to write down our findings on some examples of a quite
different nature, namely deformed dg algebras. We are very grateful
for the stimulating correspondence we had with him on the subject,
and for his interesting comments, suggestions and corrections
concerning a preliminary draft of this paper. We also thank the
referee for his quick and extremely careful reading of the
manuscript his many helpful suggestions for improvements.

\section{The homotopy category of a curved dg algebra}

\subsection{Curved dg algebras, modules and morphisms}
Curved dg algebras and modules were introduced in \cite{positselski}. We recall the definitions.
Let $k$ be a commutative ring. A cdg $k$-algebra $A$ (cdg algebra for short) consists of a graded $k$-algebra
$A = (A^i)_{i \in \Z}$,
a graded derivation $d: A \lra A$ of degree 1, and an element $c \in A^2$ with $d(c) = 0$ satisfying
$$d^2(a) = [c,a] = ca - ac$$
for all $a \in A$.
The element $c$ is called the \emph{curvature} of $A$, and $d$ is called the \emph{predifferential}. Obviously, a cdg algebra with $c = 0$ is nothing but a dg algebra.

A (left) module $M$ over a cdg algebra $A$ consists of a graded (left) $A$-module $M = (M^i)_{i \in \Z}$ endowed with a derivation $d_M: M \lra M$ of degree 1 (i.e. a degree 1 morphism with $d_M(am) = d(a)m + (-1)^{|a|}ad_M(m)$) such that
$$d_M^2(m) = cm$$
for all $m \in M$.

Modules over a cdg algebra $A$ form an abelian category $\Mod(A)$, with the obvious degree zero morphisms commuting with the predifferentials. In particular, the ground ring $k$ will be considered as a dg and cdg algebra concentrated in degree zero, and consequently $\Mod(k)$ denotes the category of complexes of ordinary $k$-modules, which we will call ``degree zero'' $k$-modules. The category of degree zero $k$-modules is denoted by $\Mod_0(k)$.

For a cdg algebra $A$, graded $A$-split exact sequences define an exact structure on $\Mod(A)$ making it into a Frobenius category. A module is projective-injective for this structure if and only if it's identity is contractible by a graded $A$-homotopy. The resulting stable category is the \emph{homotopy category} $\underline{\Mod}(A)$. Equivalently, the homotopy category $\underline{\Mod}(A)$ is obtained as the zero cohomology of the natural dg category of cdg modules.

Between cdg algebras, different kinds of morphisms can be considered. In this paper, we will only use \emph{strict} morphisms, which are a special case both of the morphisms considered in \cite{positselski}, and the morphisms of curved $A_{\infty}$-algebras considered in \cite[\S 4]{nicolas}.
A strict morphim $f: A \lra A'$ between cdg algebras is a degree zero morphism of graded algebras, commuting with the predifferentials, and preserving the curvature, i.e. with $f(c) = c'$. Cdg $k$-algebras with strict morphisms constitute a category $\mathsf{Cdg}(k)$.
A strict morphims $f: A \lra A'$ induces a restriction of scalars functor
$$\Mod(A') \lra \Mod(A).$$
Since an $A'$-homotopy can be regarded as an $A$-homotopy using $f$, we also obtain an induced restriction of scalars functor
$$\underline{\Mod}(A') \lra \underline{\Mod}(A).$$

\subsection{The initial cdg algebras}\label{parinit}
The first type of cdg algebras we consider will be called the \emph{initial} cdg algebras because each one of them is initial in a certain full subcategory of $\mathsf{Cdg}(k)$.
First we consider the cdg algebra $k[c]$ where $c$ is an element of degree $2$, the curvature (and where the predifferential is necessarily zero). This cdg algebra is clearly initial in $\mathsf{Cdg}(k)$. We also consider the cdg algebras $k[c]/c^n$ for $n > 0$. The cdg algebra $k[c]/c^n$ is initial among the cdg algebras $A$ whose curvature $c_A$ satisfies $c_A^n = 0$. Modules over $k[c]$ are precomplexes of degree zero $k$-modules, i.e. graded $k$-modules $M$ together with a predifferential $d_M: M \lra M$ satisfying no further condition. Indeed, such a precomplex $M$ can be uniquely made into a $k[c]$-module by putting $cm = d_M^2(m)$ for all $m \in M$.
Similarly, modules over $k[c]/c^n$ are precomplexes with $d^{2n} = 0$. Modules over $k[c]/c = k$ are of course ordinary complexes.

These cdg algebras are organized in the following way:
 $$k[c] \lra \dots \lra k[c]/c^n \lra k[c]/c^{n-1} \lra \dots \lra k[c]/c = k.$$
Consequently, we obtain a chain of module categories
 $$\Mod(k[c]) \longleftarrow \dots \longleftarrow \Mod(k[c]/c^n) \longleftarrow \Mod(k[c]/c^{n-1}) \longleftarrow \dots \longleftarrow \Mod(k).$$
and a chain of homotopy categories
$$\underline{\Mod}(k[c]) \longleftarrow \dots \longleftarrow \underline{\Mod}(k[c]/c^n) \longleftarrow \underline{\Mod}(k[c]/c^{n-1}) \longleftarrow \dots \longleftarrow \underline{\Mod}(k).$$
Moreover, for $A = k[c]$ or $A = k[c]/c^n$, a map of $A$-modules is contractible by a graded $A$-homotopy if and only if it is contractible by a graded $k$-homotopy (indeed, $hd + dh = f$ and $fd = df$ imply $hd^2 = d^2h$). So if we look at the chain of module categories above, the notion of contractibility is independent of the module category.
Now let $X$ be an arbitrary degree zero $k$-module. Consider the precomplexes
$$X_1 = (0 \lra X \lra 0)$$
$$X_2 = (0 \lra X \lra X \lra 0)$$
$$X_n = (0 \lra X \lra X \lra \dots \lra X \lra 0)$$
$$X_+ = (0 \lra X \lra X \lra \dots \lra X \lra X \lra \dots)$$
$$X_- = (\dots \lra X \lra X \lra \dots \lra X \lra X \ra 0)$$
$$X_{\infty} = (\dots \lra X \lra X \lra \dots \lra X \lra X \ra \dots)$$
where the maps $X \lra X$ are identities and where, for $X_n$ and $X_+$, the first non-zero entry from the left is in degree zero, and for $X_-$, the first non-zero entry from the right is in degree zero.

\begin{proposition}\label{vanishing of some precomplexes}
If $X$ is a non-zero degree zero $k$-module, then
$X_n$ is contractible if and only if $n$ is even or $n \in \{ +,  -, \infty \}$.
\end{proposition}

\begin{proof} This is a matter of alternating $0$ and $\id$ as maps $h_{i}$ in a (candidate) contracting homotopy.
\end{proof}

\section{``Derived'' categories via ``acyclic objects''}\label{parmain}

\subsection{``Derived'' categories via ``acyclic'' objects}\label{parreason}
Since cdg algebras, and modules over them, have pre\-differentials whose square is different from zero, they fail to have cohomology objects. Consequently, it is impossible to define a derived category in the usual way. In this section, we will list some possible requirements for alternative ``derived'' categories for a cdg algebra $A$.

The first, basic requirement will be that we obtain the ``derived'' category $D_{?}(A)$ as a triangle quotient of $\underline{\Mod}(A)$ by a thick subcategory $\aaa_{?}$ of ``acyclic'' objects. In fact, it seems  like this basic requirement is already largely responsible for the weird phenomena we will describe later on (see also \S \ref{parnofree}).

Next we can list possible requirements for $\aaa_{?}$, which are fulfilled in the case of the ordinary derived category of a dg algebra:

\begin{enumerate}
\item[(A1)] $\aaa_?$ contains all totalizations of short exact sequences in the abelian category $\Mod(A)$.
\item[(A2)] $\aaa_?$ is closed under coproducts.
\item[(A3)] $\aaa_?$ is closed under products.
\end{enumerate}

\begin{remark}\label{relevant consequence of A1}
Condition (A1) implies that every short exact sequence in $\Mod(A)$ gives rise to a triangle in $D_?(A)$. In fact, this is all we will need when using (A1) in this paper.
\end{remark}

Our first result is that this list of requirements makes the ``derived'' categories of all the initial cdg algebras vanish (except for $A = k$).

\begin{proposition}\label{derivedzero}
For the initial cdg algebras $A = k[c]$ or $A = k[c]/c^n$ with $n >1$, the only ``derived'' category satisfying (A1), (A2) and (A3) is $D_?(A) = 0$. If $k$ is a field, the same conclusion holds for every ``derived'' category satisfying (A1) and (A2).
\end{proposition}

\begin{proof}
If $k$ is a field, the precomplexes $k_i$ with $i \in \N_0 \cup \{+, -, \infty\}$ (see \S \ref{parinit}) that exist in $\Mod(A)$, and their shifts, are the indecomposable objects in $\Mod(A)$. By (A2), it suffices that they are acyclic. The only non finite indecomposables are contractible, hence it suffices to show that the finite indecomposables are acyclic. So by Lemma \ref{lemmaa}, in both cases, it suffices to show for $X \in \Mod_0(k)$ that $X_1$ is acyclic. Consider the exact sequence
$$0 \lra X_2[-1] \lra X_3 \oplus X_1[-1] \lra X_2 \lra 0$$ given by
$$\xymatrix{0 \ar[r] & 0 \ar[d] \ar[rr] && X \ar[d]_{\scriptsize{\left[\begin{array}{cc}\id & 0\end{array}\right]^t}} \ar[rr]^{\id} && X \ar[d]_{\id} \ar[r] & 0\\
0 \ar[r] & X \ar[d]_{\id} \ar[rr]_-{\scriptsize{\left[\begin{array}{cc}\id & \id\end{array}\right]^t}} && {X \oplus X} \ar[d]_-{\scriptsize{\left[\begin{array}{cc}\id & 0\end{array}\right]}} \ar[rr]_{\scriptsize{\left[\begin{array}{cc}\id & -\id\end{array}\right]}} && X \ar[d] \ar[r] & 0\\
0 \ar[r] & X \ar[rr]_{\id} && X \ar[rr] && 0 \ar[r] & 0}$$
From the acyclicity of $X_2$, it follows by (A1) and thickness of $\aaa_?$ that both $X_1$ and $X_3$ are acyclic. This finishes the proof.
\end{proof}

\begin{lemma}\label{lemmaa}
Let $A$ be as above. Suppose $D_?(A)$ satisfies (A1), (A2) and every object $X_1$ for $X \in \Mod_{0}(k)$ is acyclic. Then every bounded above precomplex in $\Mod(A)$ is acyclic. If $D_?(A)$ moreover satisfies (A3), then $D_?(A) = 0$.
\end{lemma}

\begin{proof}
For finite precomplexes, the proof is by induction on the length of the precomplex using (A1). Every bounded above (resp. below) precomplex can be written in $\underline{\Mod}(A)$ as a cone of coproducts (resp. products) of finite precomplexes. Using (A1) we also get the unbounded precomplexes.
\end{proof}

\begin{remark}
Note that the proof of Proposition \ref{derivedzero} makes use of the existence of $X_3$ in all of the categories $\Mod(A)$ considered. For $A = k[c]/c = k$, the classical derived category $D(k)$ is a non zero ``derived'' category satisfying (A1), (A2), (A3) (corresponding to the fact that $X_3$ does not exist and $X_1$ is not acyclic).
\end{remark}

\subsection{``Derived'' categories and base change}\label{parbase}
Another type of requirement involves the behaviour of ``acyclic'' objects, and hence ``derived'' categories, under base change. Consider a strict morphism $f: A' \lra A$ of cdg algebras, and the induced restriction of scalars functor
$$f^{\ast}: \underline{\Mod}(A) \lra \underline{\Mod}(A').$$
We can now formulate a weak and a strong base change property:
\begin{enumerate}
\item[(Bw)] The functor $f^{\ast}$ preserves ``acyclic'' objects, i.e $f^{\ast}(\aaa_?) \subseteq \aaa'_?$.
\item[(Bs)] The functor $f^{\ast}$ preserves and reflects ``acyclic'' objects, i.e. $\aaa_? = {f^{\ast}}^{-1}(\aaa'_?)$.
\end{enumerate}

Clearly, as soon as (Bw) holds, we obtain an induced restriction of scalars functor
$$f^{\ast}: D_?(A) \lra D_?(A'),$$
and if moreover (Bs) holds, this functor reflects isomorphisms.

Our next observation is that the strong base change condition combined with the conditions of \S \ref{parreason} makes all ``derived'' categories vanish.

\begin{proposition}
Let $A$ be a cdg algebra with ``derived'' category $D_?(A)$, and suppose the unique morphism $f: k[c] \lra A$ satisfies (Bs) with respect to $D_?(A)$ and $D_?(k[c]) = 0$. Then $D_?(A) = 0$.
\end{proposition}

\begin{proof}
This is obvious.
\end{proof}

\begin{example}
Consider the canonical $A \lra k$ for $A$ as in Proposition \ref{derivedzero}. Then this morphism does not satisfy (Bs) with respect to the usual derived category $D(k)$ and the ``derived'' category $D_?(A) = 0$ since $k_1$ is not acyclic in $\Mod(k)$ but becomes ``acyclic'' in $\Mod(A)$.
\end{example}

\subsection{``Derived'' categories of deformations} \label{parzeroder}
An important source of cdg algebras is given by deformations of dg algebras. Let $(A, m_{A}, d_{A})$ be a dg $k$-algebra. The Hochschild complex $\CC(A)$ is the product total complex of the double complex with
$$\CC^{\ast, n}(A) = \Hom^{\ast}_k(A^{\otimes n}, A)$$
and the familiar Hochschild differential. Consequently, a Hochschild $2$-cocycle $\phi = (\phi_n)_{n \geq 0}$ is determined by elements $\phi_0 \in A^2$, $\phi_1: A \lra A$ of degree $1$, $\phi_2: A \otimes_k A \lra A$ of degree $0$ and so on. If we concentrate on a cocycle $\phi = (\phi_0, \phi_1, \phi_2)$, this determines a first order deformation $A_{\phi}[\epsilon]$ of $A$ which is a cdg $k[\epsilon]$-algebra with multiplication $m_A + \phi_2\epsilon$, predifferential $d_A + \phi_1\epsilon$, and curvature $\phi_0\epsilon$ (a general cocycle determines a curved $A_{\infty}$-deformation, see \cite{lowen7} and \cite{kellerlowen}).

The deformation theory of algebras (\cite{gerstenhaber, gerstenhaber1}) and of abelian categories (\cite{lowenvandenbergh1, lowenvandenbergh2}) suggests that deformation should somehow take place on the derived level.

We thus wonder whether there exists a non zero ``derived'' category of $A_{\phi}[\epsilon]$ which satisfies (A1), (A2) and perhaps (A3).
First of all, note that the argument of Proposition \ref{derivedzero} for the contrary fails. Indeed, since the curvature of $A_{\phi}[\epsilon]$ is $c = \phi_{0}\epsilon$, $d_M^2$ of an $A_{\phi}[\epsilon]$-module $M$ has to factor through $\epsilon M$ so $X_3$-type objects can never exist.

Secondly, the question we ask is not complete, for we are not looking for an arbitrary derived category of $A_{\phi}[\epsilon]$, but for one that ``deforms'' $D(A)$ in some sense. A basic requirement in this respect seems to be that the strict morphism $A_{\phi}[\epsilon] \lra A$ satisfies the base-change property (Bs) with respect to $D_?(A_{\phi}[\epsilon])$ and the usual derived category $D(A)$. If this requirement is fulfulled, we say that $D_?(A_{\phi}[\epsilon])$ \emph{deforms} $D(A)$.

In the remainder of this chapter we discuss two examples where such a derived deformation does not exist.

\subsection{The cdg algebras $R_{\rho}[u]$ and $R_{\rho}[u, u^{-1}]$}\label{gradedfield}
We now introduce the two types of cdg algebras we will use. Let $R$ be a (degree zero) $k$-algebra and let $\rho \in R$ be a central element. Then $R_{\rho}[u]$ is the cdg algebra $$R[u] = (0 \lra R \lra 0 \lra Ru \lra 0 \lra Ru^2 \lra \dots)$$ where $u$ is a variable of degree $2$, with curvature $c = \rho u$.
Modules over $R_{\rho}[u]$ are precomplexes $M$ of $R$-modules with a distinquished map of precomplexes $u_M: M \lra M[2]$ for which $d_M^2 = \rho u_M$. Maps $f: (M, u_M) \lra (N, u_N)$ have to satisfy  $u_N f = f u_M$.

The localisation $R_{\rho}[u, u^{-1}]$ of $R_{\rho}[u]$ is the cdg algebra $$R[u, u^{-1}] = (\dots \lra Ru^{-1} \lra 0 \lra R \lra 0 \lra Ru \lra \dots)$$ with curvature $c = \rho u$. Modules over $R_{\rho}[u, u^{-1}]$ are modules over $R_{\rho}[u]$ where $u_M: M \lra M[2]$ is an isomorphism of precomplexes.
Up to isomorphism, they are given by precomplexes
$$\xymatrix{{\dots} \ar[r]  & M \ar[r]_{d_0} & N \ar[r]_{d_1} & M \ar[r]_{d_0} & N \ar[r] & {\dots}}$$
with $d_1 d_0 =\rho_M$ and $d_0 d_1 =\rho_N$.

We put $R[u] = R_0[u]$ and $R[u, u^{-1}] = R_0[u, u^{-1}]$.

\subsection{Some ``derived'' categories of deformations}
Consider the following two examples of $k[\epsilon]$-deformations in the diagram on the right:

$$\xymatrix{{k[\epsilon]} \ar[d] & {A_u[\epsilon] = (A[\epsilon], c = u\epsilon)} \ar[d] & {k[\epsilon]_{\epsilon}[u]} \ar[d] \ar[r] & {k[\epsilon]_{\epsilon}[u, u^{-1}] \ar[d]}\\
k & {A = (A, c = 0)} & {k[u]} \ar[r] & {k[u, u^{-1}]}}$$

In \cite[Proposition 3.13, Example 3.14]{kellerlowen}, it was shown that the ``graded field'' $k[u, u^{-1}]$ has no Morita deformation corresponding to the Hochschild cocycle $\phi = u$. Our next proposition shows that in fact, it has no reasonable corresponding ``derived'' deformation either.

\begin{proposition}\label{prophor}
For $A = k[u]$ or $A = k[u, u^{-1}]$, there is no ``derived'' category of $A_u[\epsilon]$ satisfying (A1) and deforming the classical derived category $D(A)$.
Moreover, if $k$ is a field, the only ``derived'' category of $k[\epsilon]_{\epsilon}[u, u^{-1}]$ satisfying (A1) and (A2) is $D_?(k[\epsilon]_{\epsilon}[u, u^{-1}]) = 0$.
\end{proposition}

\begin{proof}
Put $B = A_u[\epsilon]$ in either case. The proof will only make use of $k[\epsilon]_{\epsilon}[u, u^{-1}]$-modules, which are  considered as $k[\epsilon]_{\epsilon}[u]$-modules in case $A = k[u]$. Suppose we have a $D_?(B)$ satisfying (A1), (A2).
Consider the exact sequence of $B$-modules
$$\xymatrix{0 \ar[r] & {M'} \ar[r]_{\varphi} & M \ar[r] & {M''} \ar[r] & 0}$$ given by (from top to bottom):

$$\xymatrix{& 0 \ar[d] & 0 \ar[d] & 0 \ar[d] \\ {\dots} \ar[r] & 0 \ar[d] \ar[r] & k \ar[d]^{\epsilon} \ar[r] & 0 \ar[d] \ar[r] & {\dots}\\
{\dots} \ar[r] & k \ar[d]_{\id} \ar[r]_{\epsilon} & {k[\epsilon]} \ar[d]^{\id} \ar[r]_{\id} & k \ar[d]^{\id} \ar[r] & {\dots}\\
{\dots} \ar[r] & k \ar[r]_0 \ar[d] & k \ar[d] \ar[r]_{\id} & k \ar[d] \ar[r] & {\dots}\\
& 0 & 0 & 0}$$
The module $M''$ is contractible hence ``acyclic''. By (A1), the sequence determines a triangle in $D_?(B)$, so $\varphi$ becomes an isomorphism in $D_?(B)$.  Now the standard $\underline{\Mod}(B)$-triangle constructed on $\varphi$ also determines a triangle in $D_?(B)$, so  the object $\mathsf{cone}(\varphi)$ is acyclic. Now $\mathsf{cone}(\varphi)$ is given by
$$\xymatrix{{\dots} \ar[r] & {k \oplus k} \ar[rr]_-{\scriptsize{\left[\begin{array}{cc}\epsilon &\epsilon\end{array}\right]}} && k[\epsilon] \ar[rr]_-{\scriptsize{\left[\begin{array}{cc}\id & 0\end{array}\right]^t}} && {k \oplus k} \ar[r] & {\dots}}$$
which is readily seen to be isomorphic to the direct sum $M'[1] \oplus M$. It follows that both $M'$ and $M$ are acyclic. The fact that $M'$ is acyclic shows that $D_?(B)$ does not deform $D(A)$. Moreover, if $k$ is a field and $A = k[u, u^{-1}]$, then by Lemma \ref{lemindec} it shows that every indecomposable $A$-module, and hence, by (A2), every $A$-module, is acyclic. But since every $B$-module can be written as an extension of $A$-modules, this finishes the proof that $D_?(B) = 0$.
\end{proof}

\begin{lemma}\label{lemindec}
Let $k$ be a field. The indecomposable objects in $\Mod(k[u, u^{-1}])$ are given by (shifts of) $$\cdots \lra k \lra 0 \lra k \lra \cdots$$ and $$\xymatrix{{\dots} \ar[r] & k \ar[r]_{\id} & k \ar[r]_0 & k \ar[r] & {\dots}}.$$ Every object decomposes as a direct sum of these.
\end{lemma}

\begin{proof}
This easily follows from some base changes.
\end{proof}

\subsection{The link with $\Z_2$-graded cdg algebras}\label{partwee}
Instead of working with $\Z$-graded cdg algebras and modules, one can consider the parallel $\Z_2$-graded theory. We will call the corresponding objects cdg$_2$ algebras and modules, and for a cdg$_2$ algebra $A$ the related module categories are denoted by $\Mod_2(A)$, $\underline{\Mod}_2(A)$, $D_{2?}(A)$.

Any $k$-algebra $R$ with given central element $\rho \in R$ yields a cdg$_2$ algebra $R_{\rho} = R \lra 0 \lra R$ with curvature $\rho$. Modules over $R_{\rho}$ are $\Z_2$-precomplexes of $R$-modules
$$\xymatrix{M \ar[r]_{d_0} & N \ar[r]_{d_1} & M}$$
with $d_1 d_0 = \rho_M$ and $d_0 d_1 = \rho_N$.

We have the following tautology:

\begin{proposition}
Let $R$ be a $k$ algebra with central element $\rho \in R$. We have a diagram
$$\xymatrix{{\Mod_2(R_{\rho})} \ar[r]^-{\sim} \ar[d] & {\Mod(R_{\rho}[u, u^{-1}])} \ar[d]\\
{\underline{\Mod}_2(R_{\rho})} \ar[r]_-{\sim} & { \underline{\Mod}(R_{\rho}[u, u^{-1}])}  }$$
in which the first line is an equivalence and the second line is a triangle equivalence.
\end{proposition}

\begin{corollary}\label{corz2}
Let $k$ be a field. The only ``derived'' category of $k[\epsilon]_{\epsilon}$ which satisfies (A1) is $D_{2?}(k[\epsilon]_{\epsilon}) = 0$.
\end{corollary}

\begin{proof}
This is just a reformulation of Proposition \ref{prophor}.
\end{proof}

\section{``Derived'' categories via ``homotopical projectives''}\label{parmainhop}

\subsection{``Derived'' categories via ``homotopical projectives''}\label{parhop}
Let $A$ be a cdg algebra and $\aaa_{?} \subseteq \underline{\Mod}(A)$ a thick subcategory. Then by localization theory, the triangle quotient $D_{?}(A) = \underline{\Mod}(A)/\aaa_{?}$ is equivalent to the full subcategory $\ppp_? \subseteq \underline{\Mod}(A)$ of $\aaa_?$-\emph{homotopical projectives}, i.e. objects $P$ with ${\underline{\Mod}(A)}(P,X) = 0$ for all $X \in \aaa_?$.

Obviously, one can go the other way round and propose a generating class $\mmm \subseteq \underline{\Mod}(A)$ of ``homotopical projectives'', and define $X \in \underline{\Mod}(A)$ to be $\mmm$-\emph{acyclic} if ${\underline{\Mod}(A)}(M[i], X) = 0$ for all $M \in \mmm$ and $i \in \Z$.

\begin{remark}\label{gradhom}
The $\mmm$-acyclic objects can be understood in a cohomological manner. For cdg $A$-modules $M$ and $N$, consider the \emph{complex} $C_M(N) = \Hom_{\mathsf{Gr}(A)}(M,N)$ of graded $A$-module maps. Its cohomology is given by $$H^i_M(N) = H^i\Hom_{\mathsf{Gr}(A)}(M,N) = \underline{\Mod}(A)(M[-i],N)$$
Consequently, $N$ is $\mmm$-acyclic if and only if $C_M(N)$ is acyclic for every $M \in \mmm$ if and only if $H^i_M(N) = 0$ for every $M \in \mmm$ and $i \in \N$.
\end{remark}

\begin{definition}
An object $M$ of $\Mod(A)$ is {\em graded small} if the covariant functor $\Hom_{\mathsf{Gr}(A)}(M,?):\mathsf{Gr}(A)\ra\Mod(k)$ preserves arbitrary coproducts.
\end{definition}

\begin{proposition}\label{mmmacyc}
Suppose $\mmm$ is a class of objects of $\Mod(A)$ that are graded projective over $A$. Then the $\mmm$-acyclic objects form a thick subcategory $\aaa_{\mmm}$ of $\underline{\Mod}(A)$ (with corresponding $D_{\mmm}(A)$)  which satisfies (A1) and (A3). If the objects of $\mmm$ are moreover graded small, then $\aaa_{\mmm}$ also satisfies (A2).
\end{proposition}

\begin{proof}
$\aaa_{\mmm}$ is triangulated since $\underline{\Mod}(A)(M,-)$ is homological. The remainder of the claim follows from Remark \ref{gradhom}.
\end{proof}

\subsection{In the absence of free modules}\label{parnofree}
For a dg algebra $A$, the classical derived category $D(A)$ is generated by the free module $A \in \underline{\Mod}(A)$. However, for a general cdg algebra $A$, there is no natural way to make $A$ itself into an $A$-module. It seems that this fact is largely responsible for the vanishing of some ``derived'' categories discussed earlier on: in general, $\underline{\Mod}(A)$ simply contains too few modules, or, more correctly, not the right kind of modules. A related observation was made in \cite[Remark 3.18]{lowen7}.

\begin{remark}
Let $R$ be a $k$-algebra and $\rho \in R$ a central element. We consider the cdg$_2$ algebra $R_{\rho}$ of \S \ref{partwee}. Let $\ppp \subseteq \underline{\Mod}(R_{\rho})$ be the class of $R_{\rho}$-modules $M \lra N \lra M$ with $M$ and $N$ projective over $R$. Sometimes, the category $\ppp$ is considered as \emph{the} derived category of $R_{\rho}$ (for instance for $R = k[x]$, see \cite{kapustinli}, \cite{kapustinorlov}, \cite{orlov}). The fact that this is a ``good'' definition in this case is due to the fact that $k[x]$ has finite global dimension (see also \S \ref{parsecond}). In general, we know from the dg case that homotopical projectivity can not be defined on the graded level, and we have seen in Corollary \ref{corz2} that one may end up with nothing at all.
\end{remark}

Proposition \ref{mmmacyc} suggests a way of obtaining ``exotic'' derived categories by replacing the (no longer existing) free module $A$ by another \emph{graded free} module. We will investigate this further in the remainder of this section.

\subsection{A cone-like construction of cdg modules}
We now describe a construction which is reminiscent of taking the cone of a map. This construction lives in the world of predifferential graded modules. A \emph{predifferential graded $k$-algebra} (pdg $k$-algebra) is a graded $k$-algebra $A$ with a derivation $d_A: A \lra A[1]$. A predifferential graded module over $A$ is a graded $A$-module $M$ with an $A$-derivation $d_M: M \lra M[1]$, the \emph{predifferential}. Morphisms are graded morphisms commuting with the predifferentials.

As usual, a map $\phi: M \lra N$ gives rise to a map $\phi[1]: M[1] \lra N[1]$ with $d_{M[1]} = -d_M$ and $\phi[1] = \phi$.

\begin{proposition}\label{pdgcone}
Let $M$ and $N$ be pdg modules over a pdg algebra $A$ and let $\phi: M \lra N[1]$ and $\varphi: N \lra M[1]$ be pdg maps. There is a pdg module $\mathsf{cone}(\phi, \varphi)$ given by $N \oplus M$ as a graded module with predifferential $$d =
\begin{pmatrix}d_N & \phi \\ \varphi & d_M \end{pmatrix}$$
The predifferential $d$ satisfies
$$d^2 = \begin{pmatrix} d_N^2 + \phi \varphi & 0\\ 0 & \varphi \phi + d_M^2 \end{pmatrix}$$
\end{proposition}

\begin{proof}
To see that $d$ is an $A$-derivation, we consider $\mu_N: A \otimes N \lra N$ and $\mu_{M}: A \otimes M \lra M$ and we compute
$$d\begin{pmatrix} \mu_N & 0 \\ 0 & \mu_{M} \end{pmatrix} = \begin{pmatrix} d_N\mu_N  & \phi \mu_{M} \\ \varphi \mu_N &  d_{M}\mu_{M}\end{pmatrix} = \begin{pmatrix} \mu_N & 0 \\ 0 & \mu_{M} \end{pmatrix} (d_A \otimes \id_{N \oplus M} + \id_A \otimes d)$$
Of course we have $$d^2 = \begin{pmatrix} d_N^2 + \phi \varphi & d_N \phi + \phi d_M\\ \varphi d_N + d_M \varphi & \varphi \phi + d_M^2\end{pmatrix} = \begin{pmatrix} d_N^2 + \phi \varphi & 0\\ 0 & \varphi \phi + d_M^2 \end{pmatrix}$$
since $\phi$ and $\varphi$ are pdg maps.
\end{proof}

For a cdg algebra $A$, the category $\Mod(A)$ of cdg $A$-modules is clearly a full subcategory of the category of pdg $A$-modules.  For every pdg $A$-module $M$, the curvature $c$ defines a map of pdg $A$-modules
$$c_M: M \lra M[2]: m \longmapsto cm$$
\begin{proposition}\label{coney}
Let $M$ and $N$ be pdg $A$-modules over a cdg algebra $A$ and let $\phi: M \lra N[1]$ and $\varphi: N \lra M[1]$ be pdg $A$-module maps. If we have
$$d_N^2 + \phi \varphi = c_N \hspace{2cm} d_M^2 + \varphi \phi = c_M$$
then $\mathsf{cone}(\phi, \varphi)$ is a cdg $A$-module.
\end{proposition}

\begin{proof}
Immediate from Proposition \ref{pdgcone}
\end{proof}

\subsection{Derived categories constructed from $A$-splittings}
We can use Proposition \ref{coney} to construct cdg $A$-module structures on graded free $A$-modules in the following way. A cocycle $\phi \in A^i$ will be identified with any corresponding map $A[j] \lra A[j+i]$ depending on the context.

\begin{definition}
Let $A$ be a cdg algebra with curvature $c \in A^2$. A \emph{splitting} for $A$ (or $A$-splitting) consists of two cocycles $\phi \in A^{1-i}$ and $\varphi \in A^{1+i}$
with $$c - d_A^2 =  \varphi \phi = \phi \varphi.$$
The cdg $A$-module $A_{\phi,\varphi}$ is by definition $\mathsf{cone}(\phi, \varphi)$ where we consider
$$ \phi: A[i] \lra A[1]  \hspace{2cm} \varphi: A \lra A[i][1]$$
\end{definition}

Since $A_{\phi, \varphi}$ is graded projective and small, we obtain a ``derived'' category $D_{\phi,\varphi}(A)$ satisfying (A1), (A2), (A3) by taking $\mmm = \{A_{\phi, \varphi}\}$ in \S \ref{parhop}.

\begin{example}
Let $A$ be an initial cdg algebra $k[c]$ or $k[c]/c^n$ for $n > 1$. Up to isomorphism, the only $A$-splitting is given by $\phi=\id$ and $\varphi = c$.  The module $A_{1,c}$ is isomorphic to $k_+= 0 \lra k \lra k \lra k \lra \dots$, which is contractible. Hence, as we already know by Propositions \ref{derivedzero} and \ref{mmmacyc}, $D_{1,c}(A) = 0$. For $A = k$, there is another $c$-splitting given by $\phi = \varphi = 0$. Here $k_{0,0} = k \oplus k[-1]$, and $D_{0,0}(k)$ is the ordinary derived category. More generally, for a dg algebra $A$, $D_{0,0}(A)$ is the ordinary derived category, whereas other $0$-splittings will yield other ``exotic'' derived categories.
\end{example}

Let us now consider an arbitrary cdg algebra $A$ with $A$-splitting $\phi \in A^{1-i}$, $\varphi \in A^{1+i}$. We will try to understand the cohomology determined by $A_{\phi, \varphi}$ by computing the differential on $C_{A_{\phi, \varphi}}(M) = \Hom_{\mathsf{Gr}}(A)(A_{\phi, \varphi}, M)$ for an arbitrary cdg $A$-module $M$. As a graded module, $$C_{A_{\phi, \varphi}}(M) \cong M \oplus M[-i]$$ and we obtain for $m \in M^j$, $n \in M^{j-i}$:
$$d(m, n) = (d_M(m) + (-1)^j \varphi n, d_M(n) + (-1)^j \phi m)$$
This yields the following notions:
the element $(m,n)$ is a cocycle if
$$d_M(m) = (-1)^{j+1}\varphi n \hspace{2cm} d_M(n) = (-1)^{j+1}\phi m$$
and the element $(m,n)$ is a boundary if there exist $h \in M^{j-1}$, $k \in M^{j-i-1}$ with
$$m = d_M(h) + (-1)^{j+1} \varphi k \hspace{2cm} n = d_M(k) + (-1)^{j+1} \phi h$$

\begin{example}\label{ku}
Consider for a $k$-algebra $R$ with central element $\rho$ the cdg algebra $A = R_{\rho}[u]$ as defined in \S \ref{gradedfield}.
We use the $A$-splitting $\phi = \rho$, $\varphi = u$ to construct $D_{\rho, u}(A)$. The object $A_{\rho, u}$ is isomorphic to
$$\xymatrix{0 \ar[r] & R \ar[r]_{\rho} & R \ar[r]_{\id} & R \ar[r]_{\rho} & R \ar[r]_{\id} & R \ar[r]_{\rho} & {\dots}}$$
First of all, note that if $\rho$ is not invertible, then the object $A_{\rho, u}$ is not contractible. Consequently, $A_{\rho, u}$ is not $A_{\rho, u}$-acyclic, and $D_{\rho, u}(A) \neq 0$.

Let us now take $\rho = 0$, so $A$ is a dg algebra. If $M$ is a module with $u_M = 0$, then clearly $M$ is acyclic if and only if $M$ is acyclic in the classical sense. But if we consider for example $M = \dots \lra R[\epsilon] \lra R[\epsilon] \lra \dots$ with differential $\epsilon$ with $u_M = \id$, we have a cocycle $(1, \epsilon)$, but we can never have $1 = \epsilon h - \epsilon k$, so $(M, u_M = \id)$ is not acyclic with respect to the splitting $(0, u)$.
\end{example}

\begin{example}\label{kuu-}
Consider $A = R_{\rho}[u, u^{-1}]$  for $\rho \in R$ as defined in \S \ref{gradedfield}.
The object $A_{\rho, u}$ is isomorphic to
$$X = \xymatrix{\dots \ar[r] & R \ar[r]_{\rho} & R \ar[r]_{\id} & R \ar[r]_{\rho} & R \ar[r]_{\id} & R \ar[r]_{\rho} & {\dots}}$$ with $u_X = 1$. This object is contractible hence $D_{\rho, u}(A) = 0$.
\end{example}

\subsection{Deformed derived categories}

Let $A_{\phi}[\epsilon]$ be a $k[\epsilon]$-deformation of a dg $k$-algebra $A$.

\begin{proposition}\label{liftgoed}
Suppose $\mmm$ is a collection of objects in $\Mod(A_{\phi}[\epsilon])$ and put $\mmm_0 = \{ k \otimes_{k[\epsilon]} M \, | \, M \in \mmm \}$ in $\Mod(A)$. Then $D_{\mmm}(A_{\phi}[\epsilon])$ deforms $D_{\mmm_0}(A)$.
In particular, if $\mmm_0$ is a collection of homotopical projectives generating $D(A)$, the result holds with $D_{\mmm_0}(A) = D(A)$.
\end{proposition}

\begin{proof}
We have $\underline{\Mod}(A_{\phi}[\epsilon])(M, N) = \underline{\Mod}(A)(k \otimes_{k[\epsilon]} M, N)$ for $M \in \mmm$ and $N \in \Mod(A)$.
\end{proof}

We will now consider a special case of deformed cdg algebras. Let $A$ be a dg $k$-algebra and $\phi \in A^2$ a cocycle. The deformed cdg algebra $A_{\phi}[\epsilon]$ over $k[\epsilon]$ is the algebra $A[\epsilon]$ with curvature $c = \phi \epsilon$. We can construct the derived category $D_{\phi, \epsilon}(A_{\phi}[\epsilon])$ using the obvious $A$-splitting. However, this derived category has to be considered as a deformation of $D_{\phi, 0}(A)$ and not of $D(A)$!

\begin{proposition}\label{propnondef}
The derived category $D_{\phi, \epsilon}(A_{\phi}[\epsilon])$ deforms $D_{\phi, 0}(A)$.
\end{proposition}

\begin{proof}
Immediate from Proposition \ref{liftgoed}.
\end{proof}

\begin{example}\label{nonzero!}
In Examples \ref{ku} and \ref{kuu-}, we can take $R = k$ over $k$ and $\rho = 0$, and we can take $R = k[\epsilon]$ over $k[\epsilon]$ or over $k$ and $\rho = \epsilon$. It follows that both $D_{u, 0}(k[u])$ and $D_{u, \epsilon}(k[u]_u[\epsilon]) = D_{u, \epsilon}(k[\epsilon]_{\epsilon}[u])$ are non zero ``derived'' categories satisfying (A1), (A2), (A3).
\end{example}

\section{Some ``derived'' categories for arbitrary cdg algebras}\label{parlit}

In this section we take a look at some specific definitions of ``derived'' categories that are defined for arbitrary cdg algebras, which have been studied in the literature.

\subsection{The bar derived category}\label{parbar}

Let $A$ be a unital cdg algebra over a commutative ring $k$. In \cite[\S 8.2]{nicolas}, the {\em bar derived category}, $D_{\mathrm{bar}}(A)$, was defined as a natural generalization of the {\em relative derived category} of a dg algebra. One can regard $D_{\mathrm{bar}}(A)$ as the triangle quotient of the category of unital cdg $A$-modules up to homotopy, $\underline{\Mod}(A)$, by the full subcategory formed by the so-called {\em bar acyclic} modules, namely, those which are contractible when regarded as curved $A_{\infty}$-modules over $A$. Also, it is useful to consider $D_{\mathrm{bar}}(A)$ as the homotopy category of $\Mod(A)$ endowed with a structure of model category constructed with the help of the bar/cobar adjunction. Let us briefly recall here how this adjunction looks like. Let $BA$ be the bar construction associated to $A$ (see \cite[\S 4]{nicolas}), which is a counital dg $k$-coalgebra, $\Com(BA)$ the category of counital dg comodules over $BA$, and take $\tau_!
 {A}: BA\rightarrow A$ to be the composition of the map $A[1]\rightarrow A\ko a\mapsto a$, with the projection $BA\rightarrow A[1]$. Then we can define an adjoint pair of functors
\[\xymatrix{ \Mod(A)\ar@<1ex>[d]^{R_{\tau_{A}}} \\
\Com(BA)\ar@<1ex>[u]^{L_{\tau_{A}}}
}
\]
as follows:

- $L_{\tau_{A}}N$ is the {\em cobar construction} of $N$, and it is defined to be the unital graded $A$-module $(A\otimes_{k}N, m_{2}^{A}\otimes \id_{N})$ endowed with the differential
\[d_{L_{\tau_{A}}N}:=d_{A}\otimes \id_{N}+\id_{A}\otimes d_{N}+(m_{2}^{A}\otimes \id_{N})(\id_{A}\otimes\tau_{A}\otimes \id_{N})(\id_{A}\otimes\Delta_{N}),
\]
where $d_{N}$ is the codifferential of $N$, $m_{2}^{A}$ is the multiplication of $A$ and $\Delta_{N}$ is the comultiplication of $N$.

- $R_{\tau_{A}}M$ is the {\em bar construction} of $M$, and it is defined to be the counital graded $BA$-comodule $(BA\otimes_{k}M, \Delta_{BA}\otimes \id_{M})$, endowed with the codifferential
\[d_{R_{\tau_{A}}M}:=d_{BA}\otimes \id_{M}+\id_{BA}\otimes d_{M}-(\id_{BA}\otimes m_{2}^{A})(\id_{BA}\otimes\tau_{A}\otimes \id_{M})(\Delta_{BA}\otimes \id_{M}),
\]
where $d_{M}$ is the differential of $M$, $d_{BA}$ is the codifferential of $BA$ and $\Delta_{BA}$ is the comultiplication.

\begin{remark}
It was proved in \cite{nicolas} that both the bar and the cobar construction admit a more conceptual definition, being solutions of universal problems. We use this approach in the proof of Lemma \ref{compatibility} below.
\end{remark}

It turns out that a unital cdg $A$-module $M$ is bar acyclic if and only if the dg $BA$-comodule $R_{\tau_{A}}M$ is contractible, that is to say, equivalent to $0$ in the category of dg $BA$-comodules up to homotopy.

\begin{lemma}\label{A2 and A3 in the bar case}
The bar acyclic cdg $A$-modules satisfy (A2) and (A3).
\end{lemma}
\begin{proof}
Notice that $R_{\tau_{A}}$ preserves products because it has a left adjoint. On the other hand, it is straightforward to check that $R_{\tau_{A}}$ also preserves coproducts.
\end{proof}

Let $f:A\ra A'$ be a morphism of unital cdg algebras. Associated to it we have an adjoint pair
\[\xymatrix{\Mod(A')\ar@<1ex>[d]^{f^{*}} \\
\Mod(A)\ar@<1ex>[u]^{A'\otimes_{A}?}
}
\]
where $f^*$ is the {\em restriction of scalars} along $f$ and $A'\otimes_{A}?$ is the {\em extensions of scalars}. We can also consider the adjoint pair
\[\xymatrix{\Com(BA)\ar@<1ex>[d]^{B(f)_{*}} \\
\Com(BA')\ar@<1ex>[u]^{BA*_{BA'}?}
}
\]
where $B(f)_{*}$ is the {\em corestriction of scalars} along the bar construction $B(f)$ of $f$ and $BA*_{BA'}?$ is the corresponding {\em coextension of scalars}.

\begin{lemma}\label{compatibility}
The following squares are commutative up to an isomorphism of functors
\[\xymatrix{\Mod(A)\ar[d]_{A'\otimes_{A}?} & \Com(BA)\ar[l]_{L_{\tau_{A}}}\ar[d]^{B(f)_{*}} \\
\Mod(A') & \Com(BA)\ar[l]^{L_{\tau_{A'}}}
}
\hspace{1cm}
\xymatrix{\Mod(A)\ar[r]^{R_{\tau_{A}}} & \Com(BA) \\
\Mod(A')\ar[u]^{f^*}\ar[r]_{R_{\tau_{A'}}} & \Com(BA')\ar[u]_{BA*_{BA'}?}
}
\]
\end{lemma}

\begin{proof}
Here we use that the bar/cobar constructions are uniquely determined (up to isomorphism of functors) by the following isomorphisms
\[\Mod(A)(L_{\tau_{A}}N,M)\cong T_{\tau_{A}}MC(\Hom^\bullet(N,M)[-1])\cong\Com(BA)(N,R_{\tau_{A}}M)
\]
natural in $N$ and $M$, where $\Hom^\bullet(?,?)$ is the internal Hom-functor in the category of graded $k$-modules, and $T_{\tau_{A}}MC(\Hom^\bullet(N,M)[-1])$ is the {\em tangent space in $\tau_{A}$ to the set of solutions of the Maurer-Cartan equation} of $\Hom^\bullet(N,M)[-1]$ regarded as a cdg module over $\Hom^\bullet(BA,A)$, which is a cdg algebra endowed with the obvious curvature, predifferential and `convolution' product (see \cite[\S 6.3]{nicolas}). Now, it is easy to prove that we have isomorphisms
\begin{align}
\Mod(A')(A'\otimes_{A}L_{\tau_{A}}N,M')\cong\Mod(A)(L_{\tau_{A}}N,f^*M')\cong \nonumber \\
\cong T_{\tau_{A}}MC(\Hom^\bullet_{k}(N,f^*M')[-1])\cong \nonumber \\
\cong T_{\tau_{A'}}MC(\Hom^\bullet_{k}(B(f)_{*}N,M')[-1])\cong \nonumber \\
\cong \Mod(A')(L_{\tau_{A'}}(B(f)_{*}N),M') \nonumber
\end{align}
natural in $N$ and $M'$, which follows from the identity $f\tau_{A}=B(f)\tau_{A'}$.
\end{proof}

To study the behaviour of bar acyclic modules with respect to the change of rings, we need  the following result:

\begin{lemma}\label{inverting maps}
Let $A'$ be a unital cdg algebra and $M$ a unital cdg $A'$-module. Suppose $\psi:BA'\otimes_{k}M\ra BA'\otimes_{k}M$ is a morphism of graded $BA'$-comodules such that $\psi(1_{k}\otimes m)=0$ for each $m\in M$. Then for each $z\in BA'\otimes_{k}M$ there exists a natural number $n\geq 1$ such that $\psi^n(z)=0$. In particular, $\id-\psi$ is an isomorphism with inverse given by $\sum_{n\geq 0}\psi^n$.
\end{lemma}

\begin{proof}
Consider the filtration
\[0\subseteq F_{0}\subseteq\dots\subseteq F_{n}\subseteq\dots BA'\otimes_{k}M,
\]
with $F_{n}:=(k\oplus A'[1]\oplus\dots (A'[1])^{\otimes n})\otimes_{k}M\ko n\geq 0$. Let $\eta:BA'\ra k$ be the counit of the coalgebra $BA'$, and denote by $\psi_{0}$ the composition of the map $p_{M}:BA'\otimes_{k}M\ra M\ko x\otimes m\mapsto \eta(x)m$, with $\psi$. Notice that $\psi=(\id_{BA'}\otimes\psi_{0})(\Delta_{BA'}\otimes \id_{M})$ and that $\psi_{0}(1_{k}\otimes m)=0$ for all $m\in M$. This implies that $\psi(F_{n})\subseteq F_{n-1}$ for each $n\geq 0$ and, in particular, $\psi^{n+1}(F_{n})=0$.
\end{proof}

\begin{proposition}\label{preserved and reflected bar acyclics}
\begin{enumerate}
\item The functor $f^{*}:\Mod(A')\ra\Mod(A)$ preserves bar acyclic modules.
\item Assume that $k$ is a field and $A'$ (and hence $A$) has a non-zero curvature. Then $f^*:\Mod(A')\ra\Mod(A)$ reflects bar acyclicity.
\end{enumerate}
\end{proposition}

\begin{proof}
(1) That $f^{*}$ preserves bar acyclic modules follows directly from the commutativity of the second square in Lemma \ref{compatibility}.

(2) {\em Case 1: the curvature of $A'$ is not nilpotent.} By using the obvious commutative triangle
\[\xymatrix{\Mod(A')\ar[rr]^{f^*}\ar[dr] && \Mod(A)\ar[dl] \\
&\Mod(k[c])&
}
\]
and part (1) of this proposition, it suffices to prove the statement for $A=k[c]$ and $f:k[c]\ra A'$ being the unique morphism of cdg algebras.

{\em Step 1.1: Construction of a morphism of graded $k$-modules $s:A'\ra k[c]$.} We claim that for each $i\geq 1$, the map
\[f^{2i}:kc^i\ra A'^{2i}\ko rc^i\mapsto rc_{A'}^{i}
\]
is injective. Indeed, if there exist an element $r\in k\setminus{\{0\}}$ such that $rc_{A'}^{i}=0$, then $c_{A'}^{i}=r^{-1}rc_{A'}^{i}=0$, which is a contradiction. Therefore, since $k$ is a field, for each $i\geq 1$ the map $f^{2i}$ is a split injection of $k$-modules, \ie there exists a morphism of $k$-modules
\[s^{2i}: A'^{2i}\ra kc^i,
\]
such that $s^{2i}f^{2i}=\id$. By taking $s^{i}:=0$ for every $i\leq 0$ and every odd $i$, we get a morphism $s:A'\ra k[c]$ of graded $k$-modules.

{\em Step 1.2: bar acyclicity reflected.} Let $M$ be a unital cdg $A'$-module and assume there exists a morphism
\[h:R_{\tau_{k[c]}}(f^{*}M)\ra R_{\tau_{k[c]}}(f^{*}M)
\]
of graded comodules homogeneous of degree $-1$ satisfying $hd+dh=\id$, where $d$ is the codifferential of $R_{\tau_{k[c]}}(f^{*}M)$. Let $s:A'\ra k[c]$ be the morphism of graded $k$-modules constructed in step 1.1 of the proof, and let $B(s):BA'\ra B(k[c])$ be the morphism induced by $s$. Define $h'_{0}$ to be the composition
\[h'_{0}:R_{\tau_{A'}}(M)\stackrel{B(s)\otimes \id_{M}}{\longrightarrow}R_{\tau_{k[c]}}(f^{*}M)\stackrel{h}{\rightarrow}R_{\tau_{k[c]}}(f^{*}M)\stackrel{p_{M}}{\rightarrow}M,
\]
where $p_{M}: R_{\tau_{k[c]}}(f_{*}M)\ra M\ko x\otimes m\mapsto \eta(x)m$, with $\eta:B(k[c])\ra k$ being the counit of the coalgebra $B(k[c])$, and take $h':R_{\tau_{A'}}(M)\ra R_{\tau_{A'}}(M)$ to be the morphism of graded comodules defined by
\[h'=(\id_{BA'}\otimes h'_{0})(\Delta_{BA'}\otimes \id_{M}).
\]
The fact that $h'$ is compatible with the comultiplication follows from the fact that we are working over a tensor coalgebra.
Let $d'$ be the codifferential of $R_{\tau_{A'}}(M)$ and put
\[\phi:=h'd'+d'h'.
\]
Since $\phi^{-1}d'=d'\phi^{-1}$, it suffices to prove that $\phi$ is invertible. Thanks to Lemma \ref{inverting maps}, this is the case if $\phi i_{M}=i_{M}$, where $i_{M}$ is the map $M\rightarrow R_{\tau_{A'}}(M)\ko m\mapsto 1_{k}\otimes m$. The identity $(\Delta_{BA'}\otimes \id_{M})\phi=(\id_{BA'}\otimes\phi)(\Delta_{BA'}\otimes \id_{M})$ is easily checked. From this it follows the identity $\phi=(\id_{BA'}\otimes p_{M}\phi)(\Delta_{BA'}\otimes \id_{M})$, which implies that $\phi i_{M}=i_{M}$ holds whenever $p_{M}\phi i_{M}=\id_{M}$. Finally, it is straightforward to check
\[p_{M}\phi i_{M}=p_{M}h'd'i_{M}+p_{M}d'h'i_{M}=p_{M}hdi_{M}+p_{M}dhi_{M}=p_{M}i_{M}=\id_{M}.
\]

{\em Case 2: the curvature $c_{A'}$ of $A'$ is nilpotent, with $c_{A'}^{n}=0$ and $c_{A'}^{i}\neq 0$ for $1\leq i\leq n-1$.} By using the obvious commutative triangle
\[\xymatrix{\Mod(A')\ar[rr]^{f^*}\ar[dr] && \Mod(A)\ar[dl] \\
&\Mod(k[c]/c^n)&
}
\]
and part (1) of this proposition, it suffices to prove the statement for $A=k[c]/c^n$ and $f:k[c]/c^n\ra A'$ being the unique morphism of cdg algebras.

{\em Step 2.1: Construction of a morphism of graded $k$-modules $s:A'\ra k[c]/c^n$.} We claim that for each $1\leq i\leq n-1$, the map
$rc^i\mapsto rc_{A}^{i}$
is injective. Indeed, if there exist an element $r\in k\setminus{\{0\}}$ such that $rc_{A}^{i}=0$ for some $1\leq i\leq n-1$, then $c_{A}^{i}=r^{-1}rc_{A}^i=0$, which is a contradiction. Then, for each $1\leq i\leq n-1$, there exists a morphism $s^{2i}$ of $k$-modules such that $s^{2i}f^{2i}=\id$. By taking $s^{j}:=0$ for $j\neq 2i\ko 1\leq i\leq n-1$, we construct a morphism $s:A'\ra k[c]/c^n$ of graded $k$-modules.

{\em Step 2.2: bar acyclicity reflected.} Similar to step 1.2.
\end{proof}

The following result suggests that bar acyclic cdg modules are abundant. Indeed, it will help us to prove in Corollary \ref{the bar derived category vanishes} below that bar acyclic cdg modules are almost everywhere.

\begin{lemma}\label{k is always zero}
\begin{enumerate}
\item $k$ is a bar acyclic cdg module over the initial cdg algebra $k[c]$ when considered as a precomplex concentrated in degree $0$.
\item $k$ is bar acyclic cdg module over the cdg algebra $k[c]/c^n$, for any integer $n\geq 2$, when considered as a precomplex concentrated in degree $0$.
\end{enumerate}
\end{lemma}

\begin{proof}
We will prove part (1), part (2) admitting a similar proof. Let us put $A:=k[c]$. We have to give a contracting homotopy of curved $A_{\infty}$-modules over the cdg algebra $A$ (see the definition in \cite[\S 6.1]{nicolas}). For this, it suffices to take
$h_{i}=0$ for each $i\neq 2$ and to take
\[
h_{2}(c^n\otimes 1_{k}):=\begin{cases}0 &\text{ if }n\neq 1, \\ 1_{k}&\text{ if }n=1.\end{cases}
\]
We have to prove that the identity
\[\id_{p}=\sum_{r+s=p}(-1)^sm^k_{1+s}(\id^{\otimes s}\otimes h_{r})+\sum_{j+k+l=p}(-1)^{jk+l}h_{j+1+l}(\id^{\otimes j}\otimes m_{k}\otimes \id^{\otimes l})
\]
holds for any $p\geq 1$, where $\id_{p}$ is the morphism
\[\id_{p}:=\begin{cases}\id_{k}&\text{ if }p=0, \\ 0 &\text{ if }p\neq 0,\end{cases}
\]
$m^k_{i}$ is the $i$th multiplication of the cdg $A$-module $k$, and
\[m_{k}:=\begin{cases}m^{A}_{k}&\text{ if }l= 0, \\ m^{k}_{k}&\text{ if }l\neq 0.\end{cases}
\]
Notice that we have $m^k_{i}=0$ for $i\neq 2\ko m^{A}_{i}=0$ for $i\notin\{0,2\}$ and $h_{i}=0$ for $i\neq 2$. Then, we only need to consider the term $m_{1+s}^k(\id^{\otimes s}\otimes h_{r})$ with $(r,s)=(2,1)$ and the term $h_{j+1+l}(\id^{\otimes j}\otimes m_{k}\otimes \id^{\otimes l})$ with $(j,k,l)\in\{(0,0,1),(1,2,0),(0,2,1)\}$.

In the case $(j,k,l)=(0,0,1)$ (corresponding to the case $p=1$) we should have the identity
\[\id_{k}=h_{2}(m_{0}^{A}\otimes \id_{k}),
\]
which obviously holds. In the case $(r,s)=(2,1)$ and $(j,k,l)\in\{(1,2,0),(0,2,1)\}$ (corresponding to the case $p=3$) we should have the identity
\[0=-m_{2}^{k}(\id\otimes h_{2})+h_{2}(\id\otimes m_{2}^{k}-m_{2}^{A}\otimes \id),
\]
which is easily seen to hold.
\end{proof}

\begin{corollary}\label{the bar derived category vanishes}
If $k$ is a field and $A$ is a unital cdg $k$-algebra with non-vanishing curvature, then $D_{bar}(A)=0$.
\end{corollary}
\begin{proof}
We have to distinguish two cases.

{\em First case: The curvature $c_{A}$ is not nilpotent.} In this case we know that, if $f: k[c]\ra A$ is the unique morphism of cdg algebras, then $f^*:\Mod(A)\ra \Mod(k[c])$ reflects bar acyclicity (see Proposition \ref{preserved and reflected bar acyclics}). Thus, it suffices to prove that every cdg $k[c]$-module (\ie every precomplex over $k$) is bar acyclic. Notice that every cdg $k[c]$-module is a coproduct of indecomposable cdg $k[c]$-modules and the class of bar acyclic cdg $k[c]$-modules is closed under coproducts. Then, it suffices to prove that every indecomposable cdg $k[c]$-module is bar acyclic. Since $k$ is a field, the precomplexes $k_i$ with $i \in \mathbb{N}_0 \cup \{+, -, \infty\}$ (see \S 2.2) and their shifts, are the indecomposable objects in $\Mod(k[c])$ up to isomorphism. Thanks to Lemma 5.6, we know that $k_{0}$ is bar acyclic. By using mapping cone constructions, we know that any $k_{n}\ko n\geq 1$, is bar acyclic. Now, $k_{-}$ (resp. $k_{+}$) can b!
 e written in $\underline{\Mod}(k[c])$ as a cone of coproducts (resp. products) of precomplexes of the form $k_{n}\ko n\geq 0$. Since the class of bar acyclic cdg $k[c]$-modules is closed under arbitrary coproducts and products, then any bounded above or bounded below precomplex is bar acyclic. Finally, in Proposition \ref{vanishing of some precomplexes} we have already proved that $k_{\infty}$ is bar acyclic (moreover, it vanishes when regarded as an object of $\underline{\Mod}(k[c])$).

{\em Second case: $c_{A}^n=0$ and $c_{A}^{i}\neq 0$ for $1\leq i\leq n-1$.} We proceed similarly, by using this time the unique morphism of cdg algebras $f:k[c]/c^{n}\ra A$.
\end{proof}

\begin{remark}
Corollary \ref{the bar derived category vanishes} also follows from the argument indicated at the end of Remark 7.3 of \cite{positselski2}.
\end{remark}

\subsection{Derived categories of the second kind}\label{parsecond}
Let $A$ be a cdg algebra.
In \cite[\S 3.3]{positselski2}, three ``derived'' categories, called \emph{derived categories of the second kind}, are considered: the \emph{absolute derived category} $D_{\mathrm{abs}}(A)$ is the universal (``largest'', corresponding to the smallest $\aaa_?$) ``derived'' category satisfying (A1), the \emph{coderived category} $D_{\mathrm{co}}(A)$ is the universal ``derived'' category satisfying (A1) and (A2), and the \emph{contraderived category} $D_{\mathrm{ctr}}(A)$ is the universal ``derived'' category satisfying (A1) and (A3).

Proposition \ref{derivedzero} yields that for the initial cdg algebras $A = k[c]$ or $A = k[c]/c^n$ with $n > 1$ over a field $k$, we have $D_{\mathrm{co}}(A) = 0$, and Proposition \ref{prophor} yields that for $A = k[\epsilon]_{\epsilon}[u, u^{-1}] = k[u, u^{-1}]_u[\epsilon]$, $D_{\mathrm{co}}(A) = 0$.

On the other hand, as soon as a cdg algebra $A$ has a non zero ``derived'' category with the correct (Ai), it follows that the corresponding derived category of the second kind is non zero as well. A concrete example where this occurs was given in Example \ref{nonzero!}. In fact, as soon as $\underline{\Mod}(A)$ contains a graded projective (resp. graded projective and graded small) object which is not contractible, we thus conclude that $D_{\mathrm{ctr}}(A)$ is (resp. $D_{\mathrm{ctr}}(A)$ and $D_{\mathrm{co}}(A)$ are) non zero. It is easy to see that all graded projective objects are ``homotopical projective'' with respect to $\aaa_{\mathrm{ctr}}$ and all graded projective graded small objects are ``homotopical projective'' with respect to $\aaa_{\mathrm{co}}$.

The following converse is due to Positselski:

\begin{theorem}[\S 3.6, 3.7 in \cite{positselski2}]\label{posthm}
Let $A$ be a cdg algebra and let $\ppp \subseteq \underline{\Mod}(A)$ be the full subcategory of \emph{graded} projective objects.

If $A$ is \emph{graded} Artinian (i.e satisfies the descending chain condition on graded submodules), then $D_{\mathrm{ctr}}(A) \cong \ppp$.

If $A$ has finite homological dimension as a graded algebra (i.e. the abelian category $\Gr(A)$ has finite homological dimension), then the absolute derived category, the coderived category and the contraderived category coincide, and they are all equivalent to $\ppp$.
\end{theorem}

That for a cdg algebra $A$ with finite graded homological dimension and zero differential, the absolute derived category $D_{\mathrm{abs}}(A)$ can be considered to be \emph{the} derived category of $A$, follows from the following well known fact:

\begin{proposition}
Let $A$ be a graded algebra with finite graded homological dimension. Then a differential graded $A$-module is homotopically projective if and only if it is graded projective. In particular, the derived category $D(A)$ is equivalent to the full subcategory $\ppp \subseteq \underline{\Mod}(A)$ of graded projective modules.
\end{proposition}

\begin{proof}
Let $P$ be a graded projective acyclic $A$-module. Consider
$$\xymatrix{ {\dots} \ar[r]_d & P \ar[r]_d & P \ar[r]_d & C \ar[r] & 0}$$ as a projective resolution of $C = \Cokern(d)$ in the category $\Gr(A)$. Since $A$ has finite homological dimension as a graded algebra, it follows that the image of $d$ is graded projective as well, whence $P$ is contractible.
\end{proof}

For a (graded) algebra $A$ with infinite homological dimension, graded projective modules need not be homotopical projective, as the example of
$$\xymatrix{{\dots} \ar[r]_{\epsilon} & {k[\epsilon]} \ar[r]_{\epsilon} & {k[\epsilon]} \ar[r]_{\epsilon} & {\dots}}$$ over $A = k[\epsilon]$ shows.

The existence of non zero derived categories of the second kind in spite of the vanishing of those categories for the initial cdg algebra $k[c]$ corresponds to the fact that the derived categories of the second kind do not satisfy the strong base change property (Bs). On the other hand, it is easily seen that the derived categories of the second kind do satisfy (Bw).

The main application of derived categories of the second kind in \cite{positselski2} is to cdg coalgebras. More precisely, for a cdg coalgebra $C$ with cdg Cobar construction $B = \mathrm{Cob}_{\omega}(C)$ (associated to a $k$-linear section $\omega$ of $C \lra k$), the author proves a beautiful ``Koszul triality'' theorem (\cite[\S 6.7]{positselski2}) in which the coderived category of $C$-comodules, the contraderived category of $C$-contramodules, and the absolute derived category of $B$-modules are proved to be equivalent. Moreover, since $B = \mathrm{Cob}_{\omega}(C)$ has finite homological dimension as a graded algebra, by Theorem \ref{posthm}, it's three derived categories of the second kind coincide.

In \cite[\S 9.4]{positselski2}, the author proves that for cofibrant dg algebras (over a ground ring of finite homological dimension), the classical derived category and all the derived categories of the second kind coincide. He also uses this fact to argue that for general dg algebras, the classical derived category and the derived categories of the second kind have to differ, as they satisfy very different functoriality properties (the classical derived category is of course invariant under classical quasi-isomorphisms of dg algebras, and every dg algebra is quasi-isomorphic to a cofibrant one).

As far as we know, there is no natural definition of a derived category of a curved dg algebra, that coincides with the classical derived category for all ordinary dg algebras.

\def\cprime{$'$}

\providecommand{\bysame}{\leavevmode\hbox to3em{\hrulefill}\thinspace}
\providecommand{\MR}{\relax\ifhmode\unskip\space\fi MR }
\providecommand{\MRhref}[2]{%
  \href{http://www.ams.org/mathscinet-getitem?mr=#1}{#2}
}
\providecommand{\href}[2]{#2}

\end{document}